\documentclass{article}

\def\1ox{{ \Omega^1_{\scriptstyle{X}} }}
\def\2ox{{ \Omega^2_{\scriptstyle{X}} }}

\def\ok1{{ \Omega^1_K }}
\def\ok2{{ \Omega^2_K }}

\def\P{{ {\bf P} }}

\def\N{{ {\bf N} }}

\def\C{{ {\bf C} }}
\def\R{{ {\bf R} }}
\def\8{{ {\infty } }}

\def\^{{ ^{\wedge} }}

\def\Z{{ {\bf Z } }}

\title{Why everyone should know number theory}
\author{Minhyong Kim}

\begin{document}
\maketitle

Pythagoras' dictum that `all is number' is well-known. 
But my impression is that even  practicing mathematicians
are often not entirely aware of the
thoroughness with which we have developed this  very idea
over the centuries in the formation of our
world-view,  from  the evolution of
the cosmos to
the structure of elementary particles, and even in
attempts to describe
the ultimate nature of human consciousness.

Any standard book on cosmology will tell you that
the universe is a manifold with such and such structure.
What is a manifold? Nothing but a collection of
coordinate charts, i.e.,
open subsets of $\R^n$ (a collection of n-tuples
of numbers) together with rules for gluing them together.
What are the rules? They are maps, each from a subset of
one coordinate chart to a subset of another coordinate
chart. What is a map? In the modern view of set
theory, a map from $U$ to $V$ is a subset of $U\times V$
satisfying certain conditions. Thus, when $U$
and $V$ are collections of tuples of numbers, so is
a map from $U$ to $V$. Thus, we see that a manifold is
nothing but a collection of collections of numbers.
It is a rather instructive exercise to try to
express in these terms extra structures on manifolds,
such as a metric, or other objects in mathematics
whose numerical nature we tend to forget (actually for good reasons).
Going from the large to the small scale, according to
quantum mechanics, the state of the entire universe
is described by a single vector in a large Hilbert
which is a tensor product in a natural way of
many other Hilbert spaces. We can completely trace the
evolution of (all objects in) the universe by keeping track
of this vector. Such a vector is of course nothing but
an $\infty$-tuple of numbers.

The most popular models of the brain these
days views it as an ensemble of discrete state
spaces, which themselves might be a collection
of binary digits. That is, the state of the brain
is a vector in $S=\{0,1\}^N$ for some
large $N$, with rules for responding
to various external
stimuli, and it is hoped that sophisticate mental states
can be modelled as certain subsets of $S$.
The problem of describing consciousness
amounts to locating a suitable region in $S$
correlated to the conscious state.
(Of course each factor of $S$ should eventually
come with a
label allowing us to interprete the different
states.)

In sum, modern science has implemented
Pythagoras' dictum to an absurd degree, certainly
beyond anything he could have imagined. (If we recall
his inspiration, on the other hand, whereby musical notes
could be constructed via suitable combinations
of  modes of string oscillations, one could argue
that
his view
was not very different from that of modern physics.)

Now, as briefly
mentioned above, there are good reasons for forgetting that
our models of the world are so dependent on numbers.
One is that any given number or collection of numbers
is too restrictive to convey the full richness of the
object we are trying to model. The fact that we can
use many different numbers, depending on our point
of view is the genesis of the important idea of
coordinate transformations. Working with
coordinate transformations for a while then naturally
leads to the idea that it might be profitable to
work with mathematical models which themselves
are not necessarily  a priori numerical in nature.
Hence, one now works with, say, abstract manifolds
which do not have a priori coordinate charts, but
which only {\em admit} coordinate charts,
or vector spaces that do not have an a priori
basis. The notion of an abstract set and its many descendants
has undoubtedly given rise to tremendous economy of description
and conceptual clarity in modern mathematics.
(According to Dieudonne, the definition of tensors
in physics texts are remnants from the days when
the only sets used in modelling were sets of numbers.)
Nonetheless, it is clear that in
practice, even abstract objects are
conceptually convenient substitutes
 for equivalence classes of numbers.

 Given the discussion above,
one might expect all mathematicians to be number theorists,
and maybe in some sense they are. However, there is
an even more
essential reason   why most practicing mathematicians
can get by with a rather naive understanding of numbers
and might be better off doing so.
This is because of the extremely useful
geometric picture of the real and complex numbers.
Much of the time, it is perfectly reasonable
to visualize the complex numbers as a geometric
plane, and base all other constructions
upon that basic picture, oblivious to 
the fine structure of our objects,
pretty much as one can do plenty of classical
physics without worrying about the fact that
the macroscopic objects we are considering
arise from the complicated interaction
of elementary particles.
Now, my claim is that the role of a number
theorist in mathematics is
exactly analogous to the role of
a particle theorist in physics. That role being
to probe the nature of the ultimate 
constituents of the
objects  that others study from a far coarser
perspective.
Thus, in contrast to the continuum picture
of the complex plane, a number theorist
is more likely to perceive of each individual
number or groups of numbers in a discrete
fashion, and in nested hierarchies reflecting
various complexities, and even attach
a symmetry group to each individual
number. It is not that number theorists avoid
the plane model, since it is also an important
tool in  much of number theory.
It is just that the plane has a much more
grainy and elaborate shape, with many levels
of microscopic detail and structure.
Now particle physicists would  be
viewed as denizens of a baroque world of their
own making if not for the tremendous impact that
the theory of the microscopic eventually has
on  the macroscopic world, and the hope that
we will someday have a complete picture of
how to reconstruct the macroscopic world from
the structure of particles and their interactions.
But can we make similar claims for number theory,
that such detailed understanding of numbers
such as the arithmetician strives for leads
to some global macroscopic insight  (e.g.  for geometry)?

Affirmative examples are abundant. Recall for
example the various impossibility proofs for
ruler and compass constructions. There the problems
posed were evidently macroscopic in nature.
The resolution had to wait for a microscopic classification
of at least certain classes of numbers in terms
of their symmetry groups: The numbers that can be
generated by ruler and compass constructions are
all algebraic numbers with solvable symmetry groups
having elementary 2-groups as the quotients of
the composition series. Doubling the cube, for example,
would have involved constructing a number, $2^{1/3}$ with
$S_3$ symmetry group.
The impossibility of squaring the circle was eventually
understood only with Lindemann's proof that $\pi$ is
in fact a transcendental number. (Recall that a
transcendental number is a number with
an infinite symmetry group.)

This classification of numbers by 
 symmetry groups, which is at present an intractable
problem, is still just a beginning of the program to
understand how to go from the microscopic to the
macroscopic.
It is not hard to `find' all numbers with simple symmetry
groups: For example, a number has symmetry group
$Z/2$ iff it can be written in the form $a+b\sqrt{n}$
where $a,b,n$ are rational. The numbers with
symmetry group $Z/3$ or $S_3$ are exactly the ones obtained
from the rationals using field operations and cube roots.
(It is slightly subtle to distinguish between the numbers with
the two different types.) 
Among the deeper results already known are
the theorem that a number can be obtained from the rationals
through field operations and radicals iff its symmetry
group is solvable, and the beautiful theorem of
Kronecker and Weber that a number can be written as  a finite Fourier
series
$$a=\sum_j q_j \exp (2\pi i k_j/n_j)$$
with rational coefficients iff its symmetry group
is abelian. 
For example,
$$\sqrt{5}=e^{2\pi i/5}-e^{4\pi i/5}-e^{6\pi i/5}+e^{8 \pi i/5}.$$

This is a prototype of a deep
`parametrization'
theorem, which tells us exactly how to generate a
class of numbers with an abstract definition. 
The vague notion I have introduced of attaining
a microscopic understanding of numbers should be
 fleshed out by
examples of this sort, until a more systematic program
can be formulated.

We note that telling numbers apart by studying their
symmetry groups is very analogous to telling
spaces apart through discrete invariants.
That is, one may be able to use geometric
intuition to see that two spaces
are homeomorphic, but to prove that they
are {\em not} homeomorphic requires the
use of discrete invariants.
Similarly with numbers, symmetry groups
provide us with convenient labels with
which to unravel
their amorphous geometry.

Leaving it to the  experts to trace the history
of such global `applications' of number theory,
we will move on to some thoroughly
modern examples:
\vspace{3mm}

Suppose $X$ is a submanifold of $\C^2$
defined by a polynomial equation $f(x,y)=0$ with
the property that the equation has only
finitely many solutions in any ring that is
finitely generated as a $\Z$-algebra.
Then $X$ admits a hyperbolic metric.
(That is, the universal covering spaces of
$X$ is the upper-half plane.)
\vspace{3mm}

Let $X$ and $Y$ be submanifolds of $\P^n$
defined by homogeneous polynomial equations
$f_1=0,\ldots , f_m=0$ and $g_1=0, \ldots, g_k=0$ 
with integral coefficients.
Suppose there is a `sufficiently large'
prime $p$ such that
$f_1=0,\ldots , f_m=0$ and $g_1=0, \ldots, g_k=0$
 have the same number of solutions
in every finite field of characteristic $p$.

Then $X$ and $Y$ have isomorphic cohomology
groups in each degree. (Grothendieck)
\vspace{3mm}

Let $X$ and $Y$ be subspaces of $C^n$ defined
by polynomial equations $f_1=0,\ldots, f_m=0$
and $g_1=0, \ldots, g_m=0$ and assume:

(1) The same monomials appear with (different)
non-zero coefficients in $f_i$ and $g_i$
for each $i$;

(2) The coefficients appearing in $f_1,\ldots, f_m$
are algebraically independent, and the same for
the $g_i$'s.

Then $X$ and $Y$ have isomorphic cohomology groups
in every degree. (Grothendieck)

After stating this theorem, I tried to find
some  concrete non-trivial examples
with many variables. I found this hard to
do because of the paucity of theorems
on algebraic independence of numbers. It is very hard to
tell if a given collection of numbers are algebraically
independent. For example, it is not known if
$e$ and $\pi$ are algebraically independent!
Thus, after millenia of research, even a number-theorist's
understanding of the structure of
numbers is embarrassingly primitive.
\vspace{3mm}

More generally, start with the
subspace $X$ of $\C^n$ given by a
set of polynomial equations $$f_1=0,\ldots, f_m=0.$$
Let $\sigma(X)$ be the space obtained
by hitting the coefficients of the the $f_i$'s
with a field automorphism $\sigma$ of $\C$.
Then $X$ and $\sigma(X)$ have isomorphic
cohomology groups. (Grothendieck)
\vspace{3mm}

 Let $X$ and $Y$ be simply connected
submanifolds of $\P^n$
defined by homogeneous polynomial equations
$f_1=0,\ldots , f_m=0$ and $g_1=0, \ldots, g_k=0$ 
with integral coefficients.
Suppose there are two sufficiently large
primes $p$ and $l$ such that the collections
$\{f_1,\ldots , f_m\}$ and $\{g_1, \ldots, g_k\}$
are congruent modulo $p$ and $l$.
Then $X$ and $Y$ have the same higher homotopy
groups (Artin-Mazur).
\vspace{3mm}

The proof of most of these results are rather
difficult, and require the machinery of arithmetical
algebraic geometry in a heavy way. I should mention
also that there are many results whose statements
are purely in complex geometry, but whose
proofs require arithmetic sophistication,
such as Mori's remarkable theorem that
a Fano variety can be covered by rational curves,
or Batyrev's theorem that two birational
Calabi-Yau Manifolds have the same Betti numbers.
Also of interest in this connection is
Deligne and Illusie's arithmetic proof
of the existence of Hodge decomposition
for projective manifolds, a result usually
thought to require hard
analysis.  Another interesting relation to 
 analysis occurs in some recent work of Igor
Rodnianski, who studies the Schr\"{o}dinger equation
on a torus and finds that the fundamental solution
is more singular at rational times, and in general,
the degree of the singularity at a given time $t$
depends on how well $t$ can be approximated by
rationals. It is clear that the precise
structure of numbers makes its influence felt
in many different areas of mathematics.

With this smattering of examples in mind,
we will move on to more general considerations,
and spend the rest of our time discussing
 a theorem
coming from 
meta-mathematics, which could be viewed as
`explaining' the ubiquity of number theory.
 There are at least two ways
to view such results: One can be a hard-core
mathematician and treat meta-mathematics
with the same disdain that the artist reserves
for the art-critic, and declare all meta-mathematical
theorems to
be irrelevant to true mathematics.
On the other hand, one can regard meta-mathematics
as a bridge linking mathematics to philosophy,
taken in the purest sense, and thence to the
eternal and most archetypal questions of human
existence. My own view tends to
be somewhere between those two extremes, but the
example I am about to discuss provides such
a nice answer to the title of this
lecture that it is hard to resist believing
in some genuine significance it must carry. 

We start the discussion with some terminological
prerequisites.

A subset $S$ of the natural numbers $\N$ is called
{\em listable} if it can be algorithmically
generated. That is, we are requiring the existence
of an algorithm which if allowed to run
indefinitely,  will eventually have any  element of
$S$ in its output set, and no element of
$\N-S$. For example, the even numbers are clearly
listable, and any of us can immediately produce
a program which will do the trick.
The fact that the prime numbers are listable is
verified in many  introductory courses in programming:
For any given number, there is an obvious algorithm
for checking whether or not it is prime.
Then one can easily design a program which goes through
the natural numbers one by one, outputs a given number if
it is prime and moves onto the next number if it is
not. When one speaks of an algorithm in this
context,  it suffices to understand the word in
an intuitive sense, and take it on faith that it
is possible to formalize this notion in a satisfactory
way. But for the morbidly curious, I will state
that an algorithm refers to a program written in the
language whose alphabet consists of the natural
numbers and a suitable collection of variables,
and a collection of labels (for commands).
The allowable commands are one of three types:
$$ \begin{array}{ccc}
V & \leftarrow &V+1 \\
V& \leftarrow & V-1\\
\mbox{IF} & V\neq 0 &\mbox{GO TO}\ \  L
\end{array}$$
This seemingly poor programming language has been
sufficient  thus far
for implementing all processes that
we feel intuitively to be an `algorithm', and it has
been shown that the algorithms arising from
this language coincide with those coming
from a wide variety of seemingly independent reasonable
models of computation.
The proposal that this will always be the
case, so that me may with confidence 
use the above
class of programs for the definition of
an algorithm is known as Church's thesis,
and is universally accepted by computer
scientists.

 Here is a rather subtle example
of a listable set: take a polynomial $f(t,x_1,\ldots, x_n)$
with integer coefficients in $n+1$ variables. For any
natural number $t$ we plug it into the first slot,
yielding a polynomial of $n$-variables. We then ask
whether or not the equation $$f(t, x_1, \ldots, x_n)=0\ \ \ \ \ \ (*)$$
has a solution in natural numbers $(x_1, \ldots, x_n)$.
Now let $S$ be exactly the set of values $t$ for
which this equation has a solution. It is
an easy exercise to check that $S$ is a listable set
(One important thing to realize is that one
does not need an algorithm that figures out
the solvability of (*) in order to list $S$).
The sets that arise in this fashion as a parameter
set for solvable Diophantine equations is called
a Diophantine set. 
Another characterization is to say that a subset of $\N$
is Diophantine if it is the collection of
first coordinates to the set of solutions
of  a Diophantine equation. 

We have just stated that every Diophantine
set is listable. Now a remarkable theorem
of Yuri Matiyasevich says that listable sets and Diophantine
sets are in fact one and the same! This is a highly non-trivial
theorem. One needs to shows that given an algorithm
that lists a set $S$ of natural numbers, one can
write down a polynomial $f(t, x_1, \ldots, x_n)$
with the property that $t\in S$ iff $(*)$ has a
solution. We call such an equation
a Diophantine representation of the
listable set.
It is instructive to try out
some elementary examples, such as
even numbers, corresponding to the
Diophantine equation $t-2x=0$, or the
squares, corresponding to
$t-x^2=0$. The composite numbers can be
represented by the equation
$$t-(x+2)(y+2)=0$$
One important example is
the collection of numbers that are powers
of a fixed number, say 2. This set is
obviously listable, but the equation
that we might naively want to use
to express it as a Diophantine set
$t-2^x=0$ is not a polynomial equation.
It is in fact non-trivial to find a
polynomial for this set, and the
study of such exponential equations
is a key ingredient in Matiyasevich's proof.
 (In fact, what he shows is
that the set of triples $a,b,c$ such
that $c=a^b$ has a Diophantine representation,
in the sense that there is a three-parameter
family of polynomials $f(t_1, t_2, t_3, x_1, \ldots, x_n)$
for which the set of triples satisfying the
exponential relation is exactly the `solvable locus'
in the parameter space. It had been shown earlier by
Julia Robinson that this statement would imply
the general form of Matiyasevich's theorem.)
If one  considers the case where
$S$ is the listable set of prime numbers
one begins to truly appreciate the uninuitive nature of the
result under discussion. But this theorem
 turns out to have
even more surprising consequences: One can reduce
the most unlikely problems to problems about
Diophantine equations! The key to such a reduction
is the possibility of encoding many collections
of mathematical objects in an algorithmic way
into natural numbers (G\"{o}del numbering in a loose sense), 
and translating statements
about the objects into statements about
subsets of $\N$.

As an example, consider the four-color theorem,
given a computer-aided proof by Appel and Hakens.
It says that for any planar map, four colors, say red, blue,
green, yellow,
are sufficient to color all its region so that
no two regions sharing a boundary have the same color.
(The map can be {\em four-colored}.)
We reduce this problem to number theory as follows:
The essence of the argument is the fact
that we can reduce the problem to combinatorics by
encoding the map into labels for the regions together
with adjacency relations, i.e., into a planar
graph whose vertices need to be colored
in a suitable fashion.
Consider the statement $P(n)$ which says that
any map with $n$ regions can be four-colored.
This is obviously an algorithimcally decidable
assertion, because the number of (some
obvious equivalence
classes of) such maps is finite and the possible
ways of choosing one of the colors for each
region is finite. We just have to check if
at least one of the ways works. Thus,
we can construct a listable subset $S$ of $\N$
by simply listing all the numbers $n$ for which
$P(n)$ is false. Now, let the equation
$f(t, x_1, \ldots, x_n)=0$
be a Diophantine representation of $S$.
Then $S$ is empty, i.e. the four-color
theorem is true, iff the Diophantine equation
has no solution. Matiyasevich's theorem
is in principle constructive, and in the case of
the four-color theorem, one can actually write
down the Diophantine equation whose unsolvability
is equivalent to it.

Recently, a possibility emerged of
reducing one of the most famous unsolved
problems in topology, the Poincare conjecture,
to the unsolvability of a Diophantine equation,
although this reduction is not quite complete.
Recall that the Poincare conjecture
 says that an orientable, compact, connected,
simply-connected 3-manifold  is actually homeomorphic
to the 3-sphere. Since every compact 3-manifold
has a finite triangulation, we need to prove
the Poincare conjecture just for finite simplicial
complexes with simplices from dimension 0 to 
dimension 3. However, one can easily enumerate
all these in an algorithmic way
(Exercise). Denote by
$\{s_0,s_1,s_2\ldots\}$
the collection of three-dimensional
simplicial complexes enumerated thus.
Out of this, one can attempt to
enumerate algorithmically exactly the three-dimensional
manifolds that satisfy the hypothesis of
the conjecture:

(1) One first enumerates the manifolds
by computing the Euler characteristic.

(2) Then for each manifold
 one computes the homology groups using the incidence
relations suitably and enumerates just those
with $H_0=H_3 =\Z$ and $ H_1=H_2=0$.

(3) Now, there
is an  algorithm for
computing the fundamental group of a simplicial
complex as well, in the sense that one can find a 
finite presentation
for it. So it might be hoped that from the list
we have arrived at after step (2), one would be able to
cull out the three-manifolds with trivial
fundamental group. Unfortunately, the process for doing
this
is unknown at present because one cannot
algorithmically decide from a presentation whether or not
a group is trivial. Nevertheless, the experts seem to
believe that such an algorithm should exist for the
presentations arising from triangulated three-manifolds.
So let us assume for the moment that we have such
an algorithm in hand.
Then we can use it to enumerate a collection of simplicial
complexes $\{c_0, c_1, \ldots\}$ which are
exactly those
satisfying the hypothesis of the Poincare conjecture.

Now the key point is this: there is an algorithm,
due to Joachim Rubinstein,
which given any finite
 three-dimensional simplicial complex will determine
whether or not it is homeomorphic to the 3-sphere.
Thus, we can go through the set of 3-d simplicial
complexes enumerated above,
and put $n\in S$ iff $c_n$ is not homeomorphic to
the 3 sphere. We arrive thereby at a listable
set $S\subset \N$. According to Matiyasevich,
there will then be  a polynomial such that
$f(t,x_1,\ldots, x_n)=0$
has a natural number solution
 in the $x_i$'s iff $t\in S$.
Thus, provided we can find an algorithm
for determining the triviality of the fundamental
group of three-manifolds,
the Poincare conjecture can be reduced to a number-theoretic
statement
 that Diophantine equation in $n+1$ variables
$$f(t,x_1,\ldots, x_n)=0$$
has no solutions in $\N$.

Before drawing conclusions about the meaning
of all this, let us consider
finally what is perhaps
an even more bizarre
consequence of Matiyasevich's theorem:
 Consider the notion
of `formalizable theory'. This is a theory made
up of a countable alphabet and a syntax, a countable
collection of axioms and rules of inference.
In the early years of this century many leading mathematicians
and philosophers believed that `all' of mathematics
is formalizable in this sense, although we have at
present rather strong grounds to be pessimistic.
 In any case, for any such
theory, it is easy to put all the well-formed
sentences in 1-1 correspondence with $\N$, i.e.,
we can number them
$\{s_0, s_1, s_2, \ldots\}$, in an algorithmic fashion.
One can then also generate all the theorems of the theory
by applying the rules of inference in some
automatic way to the axioms. That is, the
theorems of the theory will give rise to
a listable set $S=\{ n\in \N| s_n \mbox{ is a theorem} \}$.
Corresponding to this, we again find a polynomial
$f(t,x_1,\ldots, x_n)$ such that
$t\in S$ iff $f(t,x_1,\ldots, x_n)=0$ is solvable.
Thus, a one-parameter family of Diophantine
equations captures all we would like to know
about the theorems of our theory.
To check whether the 100th sentence in our theory
is a theorem, we need `only' check if
$f(100, x_1,\ldots, x_n)=0$
has a solution is natural numbers.
As an aside, we remark that
these considerations can be used to
see that there is a Diophantine equation for
which the existence of a solution is undecidable.

What are we to make of all
this?
We leave it to the audience to decide for themselves what
Matiyasevich's theorem says about the role of
number theory in mathematics. It should be admitted
in this connection that  thus far, no significant
application has been found of this technique of `reduction to a Diophantine
equation' to other areas of mathematics.
Rather, it is usually invoked to
illustrate the difficulty of the
subject of Diophantine equations.
In the absence of more intricate
developments from the theorem
within the domain of mathematics
proper, the importance to be attached to it
will remain a matter of personal taste. 
However, in conjunction with results of the sort
stated earlier, involving  straightforward applications
of number-theoretic insight, it is hard not to
believe that there is something about number theory
that is very fundamental and all-pervasive within
the entire universe of mathematical objects.
For years, I've felt the need to deny  the popular
conception of mathematics that equates it
with the study of numbers. It is only recently
that I'm returning to a suspicion that mathematics
is perhaps about numbers after all.
It is said that in ancient Greece the
comparison of large quantities was regarded
as a very difficult problem. So it was debated
by the best thinkers of the era whether there
were more grains of sand on the beach or
more leaves on the trees of the forest.
Equipped now with  systematic notation
and fluency in the arithmetic of
large integers, it is a straightforward  (albeit tedious) matter
for even a schoolchild
 to give
an intelligent answer to such a question.
At present, our understanding of the complex numbers
is about as primitive as the understanding of
large integers was in ancient Greece.
It is tempting to speculate that a precise and
systematic knowledge of the structure
of complex numbers (in particular a good 
notational system which would reflect 
arithmetic properties)
will eventually 
render many  problems of seemingly
insurmountable conceptual difficulty
trivial to schoolchildren in the far future.

The title of this lecture is obviously a bit of
an exaggeration. It is no more true
that everyone `should'
delve deeply into number theory than that
 a solid-state physicist should  have
a comprehensive knowledge of mathematics. 
A vigorous tension between people who feel the
constant need to delve into fundamentals
and those who gaily apply intuitive reasoning
is an essential one for the health of any science.
However, I will finish with a personal
recollection of my encounter with the physicist
Victor Weisskopf, whose lecture I attended in
my first year of graduate school when I was
still contemplating a research career in mathematical physics.
I asked Weisskopf how much mathematics a physics student
needs to know, to which he answered with a smile: `More.'
The general mathematician who has ever wondered
how much knowledge of numbers is necessary for
research could do far worse than
heed the same advice.

\vspace{5mm}

{\footnotesize DEPARTMENT OF MATHEMATICS, UNIVERSITY OF ARIZONA,
TUCSON, AZ 85721}

\end{document}